\documentclass[12pt]{article}
\usepackage{amssymb}
\newcommand{\la}{\langle}
\newcommand{\ra}{\rangle}

\newcommand{\ptl}{\partial}

\newcommand{\G}{\Gamma}

\newcommand{\vs}{\varsigma}

\newcommand{\ad}{\mbox{ad}}

\newcommand{\ol}{\overline}

\def \N{\hbox{$I\hskip -4pt N$}}

\def\qed{\hfill \hfill \ifhmode\unskip\nobreak\fi\ifmmode\ifinner
         \else\hskip5pt\fi\fi
 \hbox{\hskip5pt\vrule width4pt height6pt depth1.5pt\hskip 1 pt}}
\def\qed{\hfill$\Box$}

\def\a{\alpha}
\def\b{\beta}
\def\d{\delta}

\def\G{\Gamma}
\def\l{\lambda}

\def\si{\sigma}

\def\sc{\scriptstyle}
\def\ssc{\scriptscriptstyle}
\def\F{\hbox{$I\hskip -4.5pt F$}}
\def \Z{\hbox{$Z\hskip -5.2pt Z$}}
\def\dis{\displaystyle}
\def\cl{\centerline}

\def\ol{\overline}

\def\rar{\rightarrow}
\def\Rar{\Rightarrow}

\def\Lar{\Leftarrow}

\def\Lra{\Leftrightarrow}
\def\AA{\mbox{$\cal A$}}
\def\DD{\mbox{$\cal D$}}
\def\WW{\mbox{$\cal W$}}

\def\bs{\backslash}

\def\rb{\raisebox}

\def\vs{\vspace*}
\def\vsp{{}}
\def\ra{\rangle}
\def\la{\langle}
\def\ni{\noindent}
\def\nn{\mbox{$\cal N$}}
\def\ff{\mbox{$\cal F$}}
\def\WW{\mbox{$\cal W$}}
\def\hi{\hangindent}
\def\ha{\hangafter}

\def\N{\mathbb{N}}
\def\Z{\mathbb{Z}}

\def\F{\mathbb{F}}
\begin{document}
\cl{\large\bf Structure of Algebras of Weyl Type\footnote {AMS
Subject Classification - Primary: 17B20,
17B65, 17B67, 17B68.
}}
\cl{(to appear in {\em Comm.~Algebra} {\bf32} (2004) 1051-1059.)}
\par
\centerline{\bf Yucai Su} \vs{3pt}\par \cl{Department of
Mathematics, Shanghai Jiaotong University}\cl{ Shanghai 200030,
China}\cl{E-mail: ycsu@sjtu.edu.cn} \vs{10pt}\par
\cl{and}\vs{10pt}
\par \cl{\bf Kaiming Zhao}
\vs{3pt}\par \cl{Institute of Mathematics, Academy of Mathematics
and system Sciences}\cl{Chinese Academy of Sciences, Beijing
100080, China}
\par\
\vs{5pt}
\par
\cl{\bf ABSTRACT}\par\ni In a recent paper by the authors, the
associative and the Lie algebras of Weyl type
$\AA[\DD]=\AA\otimes\F[\DD]$ were introduced, where $\AA$ is a
commutative associative algebra with an identity element over a
field $\F$ of any characteristic, and $\F[\DD]$ is the polynomial
algebra of a commutative derivation subalgebra $\DD$ of $\AA$. In
the present paper, a class of the above associative and Lie
algebras $\AA[\DD]$ with $\F$ being a field of characteristic $0$
and $\DD$ consisting of locally finite derivations of $\AA$, is
studied. The isomorphism classes of these associative and Lie
algebras are determined. The structure of these algebras is
described explicitly. \vs{4pt}\par %
\ni{\it Key Words:} Lie algebra of Weyl type, associative algebra
of Weyl type, derivation, isomorphism class.
%
%
%
%
%
\par\
\vs{5pt}
\par
\cl{\bf 1. INTRODUCTION}
\par
In a recent paper by Su and Zhao (2001), the associative and the
Lie algebras of Weyl type $\AA[\DD]=\AA\otimes\F[\DD]$ are
introduced, where $\AA$ is a commutative associative algebra with
an identity element over a field $\F$ of any characteristic, and
$\F[\DD]$ is the polynomial algebra of a commutative derivation
subalgebra $\DD$ of $\AA$. The necessary and sufficient conditions
for them to be simple are also given. Precisely, the associative
algebra or Lie algebra $\AA[\DD]$ (modulo its center, as a Lie
algebra) is simple if and only if $\AA$ is $\DD$-simple (see also
Zhao, 2002).
\par
In the present paper, we study a class of associative and Lie
algebras of the above type $\AA[\DD]$ with $\F$ being a field of
characteristic $0$ and $\DD$ consisting of locally finite
derivations of $\AA$. The isomorphism classes of these associative
and Lie algebras are determined. The structure of these algebras
is described explicitly (the special case when $\DD$ consists of
locally finite but not locally nilpotent derivations of $\AA$ was
solved in Su and Zhao, 2002). The 2-cocycles and derivations of
these algebras were determined in Su (2002, 2003), and the
derivations and automorphisms of these algebras in the spacial
case when $\AA$ is the Laurent polynomial algebra, were determined
in Zhao (1993, 1994).
\par
Throughout this paper we assume that  $\F$ is a field of
characteristic zero, all vector spaces are over $\F$. Denote by
$\Z$ the ring of integers and by $\N$ the additive semigroup of
numbers $\{0,1,2,...\}$. {}From Su et al (2000), we know that the
pairs $(\AA,\DD)$, where $\AA$ is a commutative associative
algebra with an identity element 1 over $\F$ and $\DD$ is a
nonzero finite dimensional $\F$-vector space of locally finite
commuting $\F$-derivations of $\AA$ such that $\AA$ is
$\DD$-simple, are essentially those constructed as follows.
\par
Let $\ell_1,\ell_2,\ell_3$ be three nonnegative integers such that
$\ell=\ell_1+\ell_2+\ell_3>0$.
Take any {\it nondegenerate} additive
subgroup $\G$ of $\F^{\ell_2+\ell_3}$, i.e., $\G$ contains an
$\F$-basis of $\F^{\ell_2+\ell_3}$. Elements in $\G$ will be
written as
$$
\a=(\a_{\ell_1+1},\cdots,\a_\ell)
\mbox{ \ or sometimes \ }
\a=(\a_1,\cdots,\a_\ell)\mbox{ \ with }\a_1=...=\a_{\ell_1}=0.
\eqno(1.1)$$
Set $J=\N^{\ell},J_1=\N^{\ell_1+\ell_2}\times\{0\}^{\ell_3}$.
Elements in $J$ will be written as
$$
\mu=(\mu_1,\cdots,\mu_\ell).
\eqno(1.2)$$
Elements in $J_1$ are also written as $\vec i=(i_1,i_2,\cdots,i_\ell)$
with $i_{\ell_1+\ell_2+1}=...=i_{\ell}=0.$ For $i\in\Z$, we denote
$$
i_{[p]}=(0,...,0,\rb{4pt}{\mbox{$^{^{\sc p}}_{\dis i}$}},0...,0)
\mbox{ \ for \ }1\le p\le\ell
, \eqno(1.3)$$ where $p\,$ over $i$ means that $i$ appears in the
$p$-th coordinate. Let $(\AA,\cdot)$ be the semigroup algebra
$\F[\G\times J_1]$ with the basis $\{x^{\a,\vec i}\,|\,(\a,\vec i)
\in\G\times J_1\}$ and the product
$$
x^{\a,\vec i}x^{\b,\vec j}=x^{\a+\b,\vec i+\vec j}, \ \
\forall\,(\a,\vec i),(\b,\vec j)\in\G\times J_1. \eqno(1.4)$$
Denote the identity element $x^{0,0}$ by 1. Define the derivations
$\{\ptl_p\,|\,1\le p\le\ell 
\}$ on $ \AA$ by
$$
\ptl_p(x^{\a,\vec i})=\a_p x^{\a,\vec i}+i_p x^{\a,\vec i-1_{[p]}}
\mbox{ \ for \ }1\le p\le\ell
, \eqno(1.5)$$ where we adopt the convention that if a notion is
not defined but technically appears in an expression, we always
treat it as zero; for instance, $x^{\a,-1_{[1]}}=0$ for any
$\a\in\G$. For convenience, we denote
$$%
x^\a=x^{\a,0},\,\,t^{\vec i}=x^{0,\vec i},\,\,t_p=t^{1_{[p]}}
\mbox{ \ for \ }\a\in\G,\,\vec i\in J_1,\,1\le p\le\ell.\eqno(1.6)
$$%
Set
$$
\DD={\rm span}\{\ptl_p\,|\,1\le p\le\ell 
\}. 
$$ Then we obtain the pair $(\AA,\DD)$ where $\DD$
is a nonzero finite dimensional $\F$-vector space of locally
finite commuting $\F$-derivations of $ \AA$ such that $\AA$ is
$\DD$-simple. Clearly, $\{a\in \AA\,|\,\DD(a)=0\}=\F$.
\par
Denote by $\F[\DD]$ the polynomial algebra of $\DD$. Then $\F[\DD]$ has
 a basis
$\{\ptl^\mu=\prod_{i=1}^\ell\ptl_i^{\mu_i}\,|\,\mu\in J\}$. The
$\F$-vector space
$$
\WW(\ell_1,\ell_2,\ell_3,\G)
{\sc\!}={\sc\!} \AA[\DD]{\sc\!}={\sc\!}
 \AA{\sc\!}\otimes{\sc\!}\F[ \DD]{\sc\!}={\sc\!}{\rm span}
\{x^{\a,\vec i}\ptl^\mu\,|\,(\a,\vec i,\mu){\sc\!}\in
{\sc\!} \G{\sc\!}\times {\sc\!}J_1{\sc\!}\times{\sc\!} J\},
\eqno(1.7)$$
forms an associative algebra with the following product,
$$
u\partial^\mu\cdot v\partial^\nu =\sum_{\l\in
J}({\ssc\,}^{\dis\mu}_{\dis\l}{\ssc\,})u\ptl^\l(v)
\ptl^{\mu+\nu-\l} \mbox{ \ for all \ } u,v\in  \AA,\ \mu,\nu\in J,
\eqno(1.8)$$ (when the context is clear, we shall omit the symbol
``$\cdot$'' in the product), where,
$$
({\ssc\,}^{\dis\mu}_{\dis\l}{\ssc\,})=\prod_{p=1}^{\ell}
({\ssc\,}^{\dis\mu_p}_{\dis\l_p}{\ssc\,}),\ \
\ptl^\l(v)=\Big(\prod_{p=1}^\ell\ptl^{\l_p}_p\Big)(v)
\mbox{ \ for \ }\mu,\l\in J,\,v\in{\cal A}.
\eqno(1.9)$$
This algebra $\F[\DD]$ is called an {\it algebra of Weyl type}.
The induced Lie algebra from the above associative algebra
is called a {\it Lie algebra of Weyl type}, the bracket is
$$
[u\partial^\mu, v\partial^\nu] =u\partial^\mu\cdot v\partial^\nu
-v\partial^\nu\cdot u\partial^\mu \mbox{ \ for \ }u,v\in
\AA,\,\mu,\nu\in J. \eqno(1.10)$$ Obviously $\F$ is contained in
the center of $\WW$. From  Su and Zhao (2001) and Zhao (2002), we
know that the associative algebra $\WW$ and the Lie algebra
$(\WW/\F,[\cdot,\cdot])$ are simple.
\par
Denote by $M_{m\times n}$ the set of $m\times n$ matrices with
entries in $\F$ and by $GL_m$ the group of invertible $m\times m$
matrices. Let $G_{\ell_2,\ell_3}=\{(^A_B\ ^0_C)\,|\,A\in
GL_{\ell_2}, B\in M_{\ell_2\times\ell_3},C\in GL_{\ell_3}\}$ which
is a subgroup of $GL_{\ell_2+\ell_3}$. Define an action of
$G_{\ell_2,\ell_3}$ on $\F^{\ell_2+\ell_3}$ by $g(\a)=\a g^{-1}$
for $\a\in\F^{\ell_2+\ell_3},g\in G_{\ell_2,\ell_3}.$ For any
nondegenerate additive subgroup $\G$ of $F^{\ell_2+\ell_3}$ and
$g\in G_{\ell_2,\ell_3}$, the set
$$
g(\G)=\{g(\a)\,|\,\a\in\G\},
\eqno(1.11)$$
is also a nondegenerate additive subgroup of $\F^{\ell_2+\ell_3}$. Denote
by $\Omega_{\ell_2+\ell_3}$ the set of nondegenerate additive subgroups of
$\F^{\ell_2+\ell_3}$. We have an action of $G_{\ell_2,\ell_3}$ on
$\Omega_{\ell_2+\ell_3}$ by (1.11). Define the moduli space
$$
{\cal M}_{\ell_2,\ell_3}=\Omega_{\ell_2+\ell_3}/G_{\ell_2,\ell_3},
\eqno(1.12)$$
which is the set of $G_{\ell_2,\ell_3}$-orbits in $\Omega_{\ell_2+\ell_3}$.
Then our main theorem of this paper is the following.
\par\ni%
{\bf Theorem 1.1}. {\it
Let $\WW=\WW(\ell_1,\ell_2,\ell_3,\G)$,
$\WW'=\WW(\ell'_1,\ell'_2,\ell'_3,\G')$.
 Then
$\WW\cong\WW'$ (as associative algebras or Lie algebras) if and
only if $(\ell_1,\ell_2,\ell_3)=(\ell_1',\ell_2',\ell_3')$ and
there exists an element $g\in G_{\ell_2,\ell_3}$ such that
$g(\G)=\G'$. In particular, there exists a one-to-one
correspondence between the set of isomorphic classes of the
(associative or Lie) algebras of the form (1.7) and the following
set:
$$
SW=\{(\ell_1,\ell_2,\ell_3,X)\,|\,(0,0,0)\ne(\ell_1,\ell_2,\ell_3)
\in\N^3,\;X\in {\cal M}_{\ell_2,\ell_3}\}. \eqno(1.13)$$ In other
words, the set $SW$ is the structure space of the (associative or
Lie) algebras of Weyl type in (1.7). }\par Thus by Su et al
(2000), we see that the structure space of the (associative or
Lie) algebras of Weyl type in (1.7) is the same as that of simple
Lie algebras of generalized Witt type constructed in Xu (2000).
\vs{5pt}\par\
\par
\cl{\bf 2. ABOUT $\WW(\ell_1,\ell_2,\ell_3,\G)$}
\par
We shall simply denote $\WW(\ell_1,\ell_2,\ell_3,\G)$ by $\WW$.
Denote
$$
\DD_1=\sum_{p=1}^{\ell_1}\F\ptl_p,\ \
\DD_2=\sum_{p=\ell_1+1}^{\ell_1+\ell_2}\F\ptl_p,\ \
\DD_3=\sum_{p=\ell_1+\ell_2+1}^{\ell}\F\ptl_p.
\eqno(2.1)$$
For any $\ptl=\sum_{p=1}^\ell a_p\ptl_p\in\DD$ and $\a\in\G$ written as in
(1.1), define
$$
\la\ptl,\a\ra=\a(\ptl)=\sum_{p=\ell_1+1}^{\ell}a_p\a_p.
\eqno(2.2)$$ Since $\G$ is a nondegenerate subgroup of
$\F^{\ell_2+\ell_3}$, there exists an $\F$-basis
$\a^{(\ell_1+1)},\cdots,\a^{(\ell)}\in\G$ of $\F^{\ell_2+\ell_3}$,
and the dual basis $d_{\ell_1+1},\cdots, d_\ell$ of $\DD_2+\DD_3$
with respect to (2.2) such that $\la\a^{(p)},d_q\ra=\d_{p,q}$ for
$\ell_1+1\le p,q\le\ell$. 
Set $d_p=\ptl_p$ if $1\le p\le\ell_1$. 
Then
$$
B=\{d^\mu=\prod_{p=1}^\ell d^{\mu_p}_p\,|\,\mu\in J\},
\eqno(2.3)$$
is a basis of $\F[\DD]$. We define a total order on $J$ by
$$
\begin{array}{lll}
\mu<\nu& \Lra& |\mu|<|\nu|\mbox{ or }|\mu|=|\nu| \mbox{ but
}\exists\,p\mbox{ such that}\vs{4pt}\\
&&1\le p\le\ell 
 \mbox{ and } \mu_p<\nu_p\mbox{ and
}\mu_q=\nu_q\mbox{ when }q<p,\end{array} \vsp\eqno(2.4)$$ where
the value $|\mu|=\sum_{p=1}^\ell\mu_p$ is called the {\it level}
of $\mu$. For convenience, we denote
$$
J_{11}=\N^{\ell_1}\times\{0\}^{\ell_2+\ell_3},\ \
J_{12}=\{0\}^{\ell_1}\times\N^{\ell_2}\times\{0\}^{\ell_3},\ \
J_2=\{0\}^{\ell_1}\times\N^{\ell_2+\ell_3}.
\eqno(2.5)$$
\par
Denote by $\ff$ and $\nn$ the sets of {\it ad}-locally finite and
{\it ad}-locally nilpotent elements in \WW\ respectively. Then we have
\par\ni
{\bf Lemma 2.1} {\it (a) $\AA+\DD\subset\ff\subset\AA[\DD_1]+\DD$,
(b) $\AA+\DD_1\subset\nn\subset\AA[\DD_1]$, where
$\AA[\DD_1]=\AA\otimes\F[\DD_1] ={\it span}\{x^{\a,\vec
i}d^\mu\,|\,\mu=(\mu_1,...,\mu_{\ell_1},0,...,0) \in J_{11}\}$.
}\par\ni {\it Proof.} For any $u,v
\in\AA,\ptl,\ptl'\in\DD 
$, $i\in\N$, we have
$$
\matrix{ (\ad(u+\ptl))^iv\hfill\!\!\!\!&=(\ad(\partial))^iv,
\vs{4pt}\hfill\cr
(\ad(u+\ptl))^i\ptl'\hfill\!\!\!\!&=-(\ad(\ptl))^{i-1}\ptl'(u)
.\hfill\cr}
$$
Note that $\AA[\DD]$, as an associative algebra, is generated
by $\AA$ and $\DD$, and that the action of $\DD$ on $\AA$ is locally finite.
Then it follows that
$$\AA+\DD\subset\ff,\,\,\,\ \ \AA+\DD_1\subset\nn.\eqno(2.6)$$
\par\ni
{\bf Claim 1.}~~ $\ff\subset\AA[\DD_1]+\DD$.

\par
Suppose $u\in\ff\bs(\AA[\DD_1]+\DD)$.
Decompose $u$ according to the ``level'' with respect to $\DD_2+\DD_3$:
$$
u=\sum_{i=0}^n u_i,\ \ u_i=\sum_{\mu\in J_2,|\mu|=i}a_\mu d^\mu,\
\ a_\mu\in \AA[\DD_1], \eqno(2.7)$$ where $n>0$ is the maximal
number such that $u_n\ne0$. So, if $n=1$, we must have
$u_n\notin\DD_2+\DD_3$, otherwise $u\in\AA[\DD_1]+\DD$. Write
$$
u_n=\sum_{(\a,\vec i,\nu)\in M_0}x^{\a,\vec i}\ptl^\nu
u^{(\a,\vec i,\nu)}_n
\mbox{ with } u^{(\a,\vec i,\nu)}_n\in\F[\DD_2+\DD_3],
\eqno(2.8)$$
where $M_0=\{(\a,\vec i,\nu)\in\G\times J_1\times J_{11}
\,|\,u^{(\a,\vec i,\nu)}_n\ne0\}$ is finite.
Let $\G_0=\{\a\,|\,\exists\,\vec i\in J_1,\nu\in J_{11},\,
(\a,\vec i,\nu)\in M_0\}$.
Choose a total ordering on $\G$ compatible with its group structure.
Let $\b$ be the maximal element in $\G_0$.
If $\G_0\ne\{0\}$, by
reversing the ordering if necessary, we can suppose $\b>0$.
Set $\vec j={\rm max}\{\vec i\in J_1\,|\,\exists\,\nu\in J_{11},
(\b,\vec i,\nu)\in M_0\}$, and set
$\eta={\rm max}\{\nu\in J_{11}\,|\,(\b,\vec j,\nu)\in M_0\}$.
\par\ni
{\it Case 1.} $\b=0$. \vs{-0.1pt}\par In this case we see that
$u_n\in{\rm span} \{t^{\vec i}\ptl^\nu d^\mu\,|\, (\vec
i,\nu,\mu)\in J_1\times J_{11}\times J_2\}$. Write
$$
u^{(\b,\vec j,\eta)}_n=\sum_{\mu\in J_2,|\mu|=n}c_\mu d^{\mu},\
c_\mu\in\F. \eqno(2.9)$$ Let $\l={\rm max}\{\mu\in
J_2\,|\,|\mu|=n,c_\mu\ne0\}$. Suppose $\l_k\ne0$ and $\l_i=0$ for
all $i>k$ and some $\ell_1+1\le k\le\ell$.
By induction on $s$ it is easy to see
that the ``highest'' term of $({\rm ad\ssc\,}u)^s x^{\a^{(k)}}$ is
$\l_k^sc_{\l}^sx^{\a^{(k)},s\vec
j}\ptl^{s\eta}d^{s(\l-1_{[k]})}\ne0$, where, the ``highest'' term
is the term $x^{\a,\vec i}\ptl^\mu d^\nu$ with nonzero coefficient
in $({\rm ad\ssc\,}u)^s x^{\a^{(k)}}$ such that $(\a,\vec
i,\mu,\nu)$ is the maximal quadruple with respect to the
lexicographical order in $\G\times J_1\times J_{11}\times J_2$.
Thus
$$
{\rm dim}({\rm span}\{({\rm ad\ssc\,}u)^sx^{\a^{(k)}}\,|\,s\in\N\})=\infty.
\eqno(2.10)$$
This contradicts the fact that $u\in\ff$. So Case 1 does not occur.
\par\ni
{\it Case 2.} $\b\ne0$.
\par
Choose elements in $\G$: $\b^{(\ell_1+1)},\cdots,\b^{(\ell)}=\b$ which form
a basis of $\F^{\ell_2+\ell_3}$ ,
and choose a basis of $\DD_2+\DD_3$:
$ \{\ptl'_{\ell_1+1},\cdots,\ptl'_\ell\}$ which is a dual basis to
$\b^{(\ell_1+1)},\cdots,\b^{(\ell)}=\b$. Write
$u^{(\b,\vec j,\eta)}_n=\sum_{\mu\in J_2,|\mu|=n}c_\mu\ptl'^{\mu},\,
c_\mu\in\F$.
If for all $\mu\in\{\mu\in J_2\,|\,|\mu|=n,c_\mu\ne0\}$, we have
$\mu_\ell=0$, then we use the arguments as in Case 1 to conclude
also a contradiction. So there is a
$\mu\in\{\mu\in J_2\,|\,|\mu|=n,c_\mu\ne0\}$ with $\mu_\ell\ne0$.
Let $\l={\rm max}\{\mu\in J_2\,|\,|\mu|=n,\mu_\ell\ne0,c_\mu\ne0\}$.
It is clear that the highest term of $({\rm ad\ssc\,}u)^sx^{2\b}$ is
$$
\b_\ell^sc_{\l}^s(2\l_\ell)(2\l_\ell+1)\cdots(2\l_\ell+s-1)
x^{(s+2)\b,s\vec j}\ptl^{s\eta}\ptl'{}^{s(\l-1_{[\ell]})}\ne0.
\eqno(2.11)$$
Thus
$$
{\rm dim}({\rm span}\{({\rm
ad\ssc\,}u)^sx^{2\b}\,|\,s\in\N\})=\infty. \eqno(2.12)$$ This
contradicts the fact that $u\in\ff$. So Case 2 does not occur
either. Therefore $\ff\subset\AA[\DD_1]+D$.
 Thus Claim 1 and (a) follow.
\par
Now we prove the second part of (b).
{}From (a) we see that $\nn\subset\AA[\DD_1]+\DD$.
For any $u+\ptl\in\AA[\DD_1]+\DD$ with $0\ne\ptl\in\DD_2+\DD_3$,
choose $\a\in\G$ such that $\a(\ptl)\ne0$. From
$({\rm ad}(u+\ptl))^s x^\a=\a(\ptl)^sx^\a\ne0$, we see that
$u+\ptl\notin\nn$. Thus $\nn\subset\AA[\DD_1]$.
 Thus (b) follows.
\qed\par
For any subset $V\subset\AA[\DD]$, we define the following two sets:
$$
\matrix{
E(V)=\{u\in\WW\,|\,
\mbox{ for any }v\in V,\,\exists\,c_v\in\F\mbox{ such that }[v,u]=c_v u\},
\vs{4pt}\hfill\cr
N(V)=\{u\in\WW\,|\,[V,u]=0\}.
\hfill\cr}
\eqno(2.13)$$
Then we have
\par\ni
{\bf Lemma 2.2}. {\it (a) $E({\cal F})=\cup_{\a\in\G}\F x^\a$, (b)
$N({\cal N})={\rm span}\{x^{\a,\vec i}\,|\,(\a,\vec i)\in\G\times
J_{12}\}$. }\par\ni%
{\it Proof}. It is easy to verify that
$$
\matrix{\dis
E(\AA+\DD)=\bigcup_{\a\in\G}\F x^\a=E(\AA[\DD_1]+\DD),
\vs{4pt}\hfill\cr
N(\AA+\DD_1)={\rm span}\{x^{\a,\vec i}\,|\,(\a,\vec i)\in\G\times J_{12}\}
=N(\AA[\DD_1]).
\hfill\cr}
$$
{}From Lemma 2.1 we deduce (a) and (b).\qed \par\
\par
\cl{\bf3. PROOF OF THE MAIN THEOREM} \vs{1pt}\par%
 Now we are ready
to prove the isomorphism theorem.
\par\ni
{\it Proof of Theorem 1.1}.
We shall use the same notation but with a prime to denote elements
associated with $\WW'$.
\par\ni
``$\Lar$'':
By assumption, $(\ell_1,\ell_2,\ell_3)=(\ell'_1,\ell'_2,\ell'_3)$
and there exists a nondegenerate
$(\ell_2+\ell_3)\times(\ell_2+\ell_3)$ matrix
$g=({}^{A\ 0}_{B\ C})\in G_{\ell_2,\ell_3}$,
where $A\in GL_{\ell_2},B\in M_{\ell_3\times\ell_2},C\in GL_{\ell_3}$,
such that
$$
\matrix{
\hfill\tau:\G\!\!\!\!&\rar\G',
\vs{4pt}\hfill\cr
\hfill\a=(\a_{\ell_1+1},\cdots,\a_\ell)\!\!\!\!&
\mapsto\a'=(\a'_{\ell_1+1},\cdots,\a'_\ell)=
(\a_{\ell_1+1},\cdots,\a_\ell)g^{-1}=\a g^{-1},
\hfill\cr}
\eqno(3.1)$$
is a group isomorphism. Set
$$
(\ol\ptl_1,...,\ol\ptl_\ell)=(\ptl'_1,...,\ptl'_\ell)
\biggl(^{\dis I_{\ell_1}\ 0}_{\dis 0\ \ g}\biggr),
\eqno(3.2)$$
where $I_{\ell_1}$ is the $\ell_1\times\ell_1$ identity matrix.
Define a linear map:
$$
\matrix{ \si:\!\!&x^\a\mapsto x'{\ssc\,}^{\tau(\a)},\ \
\ptl_p\mapsto\ol\ptl_p,\ \ \forall\,\a\in\G,\,1\le p\le\ell
, \vs{4pt}\hfill\cr& (t_1,\cdots,t_{\ell_1+\ell_2})\mapsto
(t_1',\cdots,t'_{\ell_1+\ell_2}) \biggl(^{\dis I_{\ell_1}\ \
0}_{\dis 0\ \ (A^{{\ssc\,}\rm t})^{-1}}\biggr), \hfill\cr}
\eqno(3.3)$$ (recall notations in (1.6)), where $A^{{\ssc\,}\rm
t}$ is the transpose of $A$. It is straightforward to verify that
$$
\si(\ptl_p)(\si(t_q))=\d_{p,q}=\ptl_p(t_q),\ \ \
\si(\ptl_p)(\si(x^{\a}))=\si(\ptl_p(x^{\a})), \eqno(3.4)$$ for
$1\le p\le\ell,\,1\le q\le\ell_1+\ell_2,
\,\a\in\G$. We extend $\si$ to a linear map $\si:\WW\to\WW'$ by
$$\si(x^{\a,\vec
i}\ptl^\mu)=\si(x^\a)\prod_{q=1}^{\ell_1+\ell_2}
\si(t_q)^{i_q}\prod_{p=1}^\ell\si(\ptl_p)^{\mu_p}\mbox{ \ for \
}(\a,\vec i,\mu)\in\G\times J_1\times J.$$ By (3.4), we see that
$\si$ is  an associative isomorphism $\si:\WW\cong\WW'$. Thus they
are also isomorphic as Lie algebras.
\par\ni
``$\Rar$'':
If $\WW$ and $\WW'$ are isomorphic as associative algebras,
then they must be isomorphic as Lie algebras, so we suppose they are
isomorphic as Lie algebras.
\par
The isomorphism $\si$ maps the {\it ad}-locally finite, {\it
ad}-locally nilpotent elements to the {\it ad}-locally finite,
{\it ad}-locally nilpotent elements, thus
$\si(\ff)=\ff',\si(\nn)=\nn'$. Thus also $\si(E({\cal
F}))\!=E({{\cal F}'})$, $\si(N({\cal N}))=N({{\cal N}'})$. By
Lemma 2.2(a), for $a\in\G$, there exists $\a'\in\G'$ such that
$\si(x^\a)\in\F x'^{\a'}$, i.e., there exists a bijection
$\tau:\a\mapsto\a'$ from $\G\rar\G'$ such that
$$
\si(x^\a)=c_\a x'^{\a'},\,\,\forall\,\,\a\in\G,
\eqno(3.5)$$
where $c_\a\in\F\bs\{0\}
$.
For any $0\ne\ptl\in\DD_2+\DD_3\subset\ff$, by Lemma 2.1, we know that
$\si(\ptl)\in \ff'\subset \DD'+\AA'[\DD'_1]$.
We can write
$$
\si(\ptl)=\ol\ptl+u'_{\ptl},\mbox{ \ where \ }
\ol\ptl\in\DD'_2+\DD'_3,\,u'_{\ptl}\in\AA'[\DD'_1]. \eqno(3.6)$$
Since ${\rm ad\ssc\,}\ptl|_{E({\cal F})}\ne0$, ${\rm
ad}(\ol\ptl+u'_{\ptl})|_{E({{\cal F}'})} ={\rm
ad\ssc\,}\ol\ptl|_{E({{\cal F}'})}\ne0$. This shows that
$\theta:\ptl\mapsto\ol\ptl$ is an injection $\DD_2+\DD_3\rar
\DD'_2+\DD'_3$. Similarly if $0\ne\ptl\in\DD_3$, then ${\rm
ad\ssc\,}\ptl|_{N({\cal N})}\ne0$ is {\it semi-simple} (i.e.,
there exists a basis of $N({\cal N})$ consisting of eigenvectors
of ${\rm ad\ssc\,}\ptl$), thus ${\rm
ad}(\ol\ptl+u'_{\ptl})|_{N({{\cal N}'})} ={\rm
ad\ssc\,}\ol\ptl|_{N({{\cal N}'})}\ne0$ is also semi-simple. This
shows that $\ol\ptl\in\DD'_3$, i.e., $\theta|_{{\cal
D}_3}:\DD_3\rar\DD'_3$. Thus we have
$\ell_2+\ell_3\le\ell'_2+\ell'_3,\,\ell_3\le\ell'_3$. Exchanging
positions between $\WW$ and $\WW'$ shows that
$\ell'_2+\ell'_3\le\ell_2+\ell_3,\,\ell'_3\le\ell_3$. Thus
$(\ell_2,\ell_3)=(\ell'_2,\ell'_3)$ and $\theta:\DD_2+\DD_3\rar
\DD'_2+\DD'_3$ is a linear isomorphism mapping $\DD_3$ onto
$\DD'_3$. Thus there exists $g=({}^{A\ 0}_{B\ C})\in
G_{\ell_2,\ell_3}$ such that
$$
\theta:(\ptl_{\ell_1+1},...,\ptl_\ell)\mapsto
(\ol\ptl_{\ell_1+1},...,\ol\ptl_\ell)=
(\ptl'_{\ell+1},...,\ptl'_\ell)g. \eqno(3.7)$$ For
$\ptl\in\DD_2+\DD_3,\a\in\G$, applying $\si$ to $\la\ptl,\a\ra
x^\a=[\ptl,x^\a]$, using (3.5), (3.6), we obtain $\la\ptl,\a\ra
c_\a x'^{\a'}= [\ol\ptl+u'_{\ptl},c_\a x'^{\a'}]=
\la\ol\ptl,\a'\ra c_\a x'^{\a'}$. Thus
$\la\ptl,\a\ra=\la\ol\ptl,\a'\ra$, from this and (3.7), we obtain
$\a'=\a g^{-1}.$\par%
 Thus it remains to prove that
$\ell_1=\ell'_1$. Observe that the {\it centralizer} of $E({\cal
F})$ in $\WW$ (which by definition and by Lemma 2.2(a) is
$\{u\in\WW\,|\,[u,x^\a]=0,\,\forall\,\a\in\G\}$) is equal to
$\AA[\DD_1]$, thus $\si(\AA[\DD_1])=\AA'[\DD'_1]$. We prove by
induction on $|\vec j|+|\mu|$ that
$$%
\si(x^{\a,\vec i+\vec j}\ptl^\mu)=\si(x^{\a,\vec i})
\prod_{p=1}^{\ell_1}\si(t_p)^{j_p}
\prod_{p=1}^{\ell_1}\si(\ptl_p)^{\mu_p} \mbox{ \ for \ }(\a,\vec
i,\vec j,\mu)\in\G\times J_{12}\times J_{11}\times
J_{11},\eqno(3.8)
$$%
(recall notations in (2.5)). If $|\vec j|+|\mu|=0$, there is
nothing to prove since both sides of (3.8) are $\si(x^{\a,\vec
i})$. Suppose inductively that (3.8) holds whenever $|\vec
j|+|\mu|<n$, where $n\ge1$. Now assume $|\vec j|+|\mu|=n$. Denote
by $A_{\vec j,\mu}$ the difference between the left-hand side and
the right-hand side of (3.8). Then for $q\le\ell_1$, we have,
$$%
\begin{array}{ll}
[\si(\ptl_q),A_{\vec j,\mu}]\!\!\!\!&\dis=\si([\ptl_q,x^{\a,\vec
i+\vec j}\ptl^\mu])-\si(x^{\a,\vec i})\prod_{^{\sc p=1}_{\sc p\ne
q}}^{\ell_1}\si(t_p)^{j_p}[\si(\ptl_q),\si(t_q)^{j_q}]
\prod_{p=1}^{\ell_1}\si(\ptl_p)^{\mu_p}\\
&\dis=j_q\si(x^{\a,\vec i+\vec
j-1_{[q]}}\ptl^\mu])-j_q\si(x^{\a,\vec i})\prod_{^{\sc p=1}_{\sc
p\ne
q}}^{\ell_1}\si(t_p)^{j_p}\si(t_q)^{j_q-1}\prod_{p=1}^{\ell_1}\si(\ptl_p)^{\mu_p}
\\ &=j_qA_{\vec j-1_{[q]},\mu}=0,\end{array}\eqno(3.9)
$$%
where the first and second equalities follows from the Leibniz
rule $[u,u_1u_2]=[u,u_1]u_2+u_1[u,u_2]$ for $u,u_1,u_2\in\WW$
(where the product $u_1u_2$ is given by (1.8)), and the last
equality follows from the inductive assumption. Similarly,
$[\si(t_q),A_{\vec j,\mu}]=0$. Thus
$$%
\begin{array}{ll}
A_{\vec
j,\mu}\!\!\!\!&\in\si(\{u\in\AA[\DD_1]\,\,|\,\,[u,t_q]=[u,\ptl_q]=0\mbox{
for }1\le q\le\ell_1\})\vs{4pt}\\
&=\si(\{u\in\AA[\DD_1]\,\,|\,\,[u,t^2_q]=[u,\ptl^2_q]=0\mbox{ for
}1\le q\le\ell_1\})\vs{4pt}\\
&=\{u'\in\AA'[\DD'_1]\,\,|\,\,[u',\si(t^2_q)]=[u',\si(\ptl^2_q)]=0\mbox{
for }1\le q\le\ell_1\},\end{array}\eqno(3.10)
$$%
(note that $\{u\in\AA[\DD_1]\,\,|\,\,[u,t_q]=[u,\ptl_q]=0\mbox{
for }1\le q\le\ell_1\}$ is in fact equal to $N({\cal N})$,
cf.~Lemma 2.2(b)). Since $|\vec j|+|\mu|=n>0$, say $\mu_r>0$ for
some $r\le\ell_1$ (the proof is similar if $j_r>0$). Then $A_{\vec
j+1_{[r]},\mu-1_{[r]}}$ belongs to the right-hand side of (3.10),
thus
$$%
0=[\si(\ptl_r^2),A_{\vec j+1_{[r]},\mu-1_{[r]}}]=2(j_r+1)A_{\vec
j,\mu}+(j_r+1)j_rA_{\vec j-1_{[r]},\mu-1_{[r]}}=2(j_r+1)A_{\vec
j,\mu},
$$%
where the second equality follows from an exactly similar argument
to that in (3.9) (and using (1.8) and (1.10)). Thus proves
$A_{\vec j,\mu}=0$, i.e., (3.8) holds.\par%
Observe from Lemma 2.2(b) that $N({\cal N})$ is the semigroup
algebra $\F[\G\times J_{12}]$, which is a domain ring and which is
also the center of the algebra $\AA[\DD_1]$. Denote by
$\F(\G\times J_{12})$ the rational field of $\F[\G\times J_{12}]$.
Regarding $\AA[\DD_1]$ as a Lie algebra over $\F[\G\times J_{12}]$
and extending it as a Lie algebra over the field $\F(\G\times
J_{12})$, then it is just the classical rank $\ell_1$ Lie algebra
of Weyl type over $\F(\G\times J_{12})$. Since $\si(\F[\G\times
J_{12}])=\F[\G'\times J'_{12}]$, we can regard the rational field
$\F(\G'\times J'_{12})$ the same as $\F(\G\times J_{12})$. By
(3.8), we obtain that $\si(au)=\si(a)\si(u)$ for $a\in\F[\G\times
J_{12}],\,u\in{\rm span}\{t^{\vec j}\ptl^\mu\,|\,(\vec j,\mu)\in
J_{11}\times J_{11}\}$, which shows that the classical rank
$\ell_1$ Lie algebra $\AA[\DD_1]$ of Weyl type over the field
$\F(\G\times J_{12})$ is isomorphic to the classical rank
$\ell'_1$ Lie algebra $\AA'[\DD'_1]$ of Weyl type. Thus
$\ell_1=\ell'_1$. This completes the proof of Theorem 1.1.
\qed\par\ \vs{5pt}
 \par
\cl{\bf ACKNOWLEDGMENT}%
This work is supported by a NSF grant 10171064 of China and two
grants ``Excellent Young Teacher Program'' and ``Trans-Century
Training Programme Foundation for the Talents'' from Ministry of
Education of China.
\par\ \vs{5pt}
 \par
\cl{\bf REFERENCES}%
\par\ni\hi6ex\ha1 %
Su, Y. (2002). 2-Cocycles on the Lie algebras of generalized
differential operators. {\it Comm.~Algebra} 30:763--782.
\par\ni\hi6ex\ha1 %
Su, Y. (2003). Derivations of generalized Weyl algebras. {\it
Science in China A} 46:346--354.
\par\ni\hi6ex\ha1 
Su, Y., Xu, X., Zhang, H.
(2000). Derivation-simple algebras and the structures of Lie
algebras of Witt type. {\it J.~Algebra}
233:642--662.%
\par\ni\hi6ex\ha1 
Su, Y., Zhao, K. (2001). Simple algebras of Weyl type.
{\it Science in China A
} 44:419--426.%
\par\ni\hi6ex\ha1 
Su, Y., Zhao, K. (2002). Isomorphism classes and automorphism
groups of algebras of Weyl type. {\it Science in China A
} 45:953--963.%
\par\ni\hi6ex\ha1 
Xu, X. (2000). New generalized simple Lie algebras of Cartan type
over
a field with characteristic 0. {\it J.~Algebra} 224:23--58.%
\par\ni\hi6ex\ha1 
Zhao, K. (1993). Lie algebra of derivations of algebras of
differential operators. {\it Chinese Science Bulletin}
38:793--798.%
\par\ni\hi6ex\ha1 
Zhao, K. (1994). Automorphisms of algebras of differential
operators. {\it J.~of Capital Normal University} 1:1--8.%
\par\ni\hi6ex\ha1 
Zhao, K. (2002). Simple algebras of Weyl type, II. {\it
Proc.~Amer.~Math.~Soc.} 130:1323--1332.
\end{document}